\documentstyle[12pt]{article}
\input epsf.tex
\begin{document}

\author{E. P. Volokitin and V. V. Ivanov
\thanks
{The research was supported by the Russian Foundation for
Basic Research (Grants 99--01--00574, 99--01--00575).}
}
\title{Isochronicity and Commutation of Polynomial
Vector Fields}
\date{}
\maketitle

\newtheorem{theorem}{Theorem}
\newtheorem{lemma}{Lemma}

\quad

\begin{center}
\parbox{14cm}{{\bf Abstract.}
We study a connection between the isochronicity of a center of a
polynomial vector field and the existence of a polynomial commuting
system. We demonstrate an isochronous system of degree 4
which does not commute with any polynomial system. We prove that among
the Newton polynomial systems only the Lienard and Abel systems may commute
with transversal polynomial fields. We give a full and constructive
description of centralizers
of the Abel polynomial systems. We give new examples of isochronous systems.
}
\end{center}

\vspace{.5cm}

A center of a planar dynamical system
is said to be isochronous
if all cycles near it have the same period.
The isochronicity problem
has a long history \cite{1}.
After two or three decades of oblivion,
the interest in this problem has been revived
partially due to proliferation of powerful methods of computerized research.
At this juncture, most attention is paid to polynomial vector fields
(see, for instance, \cite{2,3} and the bibliography therein).
The characteristic feature of the present-day activity
is the fact that some long forgotten results
are rediscovered and at the same time
new nontrivial examples of isochronous systems appear,
giving rise to fresh and ingenious ideas.

\begin{center}
{\large \bf A~Dynamical Meaning of the Cauchy--Riemann  Conditions}
\end{center}

One of these ideas leans upon a remarkable observation by
Villarini \cite{4} which allowed him to gave an elegant proof
to a well-known theorem by Lukashevich \cite{1,5}.
Lukashevich's theorem asserts that
each center of a system
\refstepcounter{equation}
\label{1}
$$
\begin{array} {ll}
\dot x &= p (x, y),\\
\dot y &= q (x, y)
\end{array}
\eqno{(\ref{1})}
$$
whose right-hand sides satisfy the Cauchy--Riemann conditions
$$ p_x = q_y, \quad p_y = -q_x$$
is isochronous.

From the viewpoint of dynamics, these conditions mean that
system (\ref{1}) commutes with the orthogonal system
\refstepcounter{equation}
\label{2}
$$
\begin{array}{ll}
\dot x &= \phantom{+} q (x, y),\\
\dot y &= - p (x, y).
\end{array}
\eqno{(\ref{2})}
$$
Indeed, the commutator, or Lie bracket \cite{6},
of the vector fields $ (p, q) $  and $ (q, -p) $
is the field with the components
$$
\begin{array}{ll}
q (p_x - q_y) - p (p_y + q_x), \\
p (p_x - q_y) + q (p_y + q_x)
\end{array}
$$
whose vanishing at a noncritical point with $ p^2 + q^2 \ne 0$
is equivalent to the fulfillment of the Cauchy--Riemann conditions.

Thus, if a point, moving according to equations (\ref{1}),
returns to the initial position at time $ T$ and
the orthogonal flow (\ref{2}) afterwards carries it
to a ~nearby orbit of (\ref{1}),
then the fact that the flows commute implies that
the motion of the point along the latter
by time $T$ also terminates at the same place where it starts.

This argument proves Lukashevich's theorem and also allows us
to make some interesting conclusions
about the properties of the system under consideration.
For example, the orbits of such a~system, lying nearby
a~ cycle, are necessarily closed, so that there are no limit cycles.
Moreover, if an analytical
system is defined over a simply connected domain then each cycle
of the system surrounds a unique singular point, which is
an isochronous center, and lies in the basin of this center
together with the region that it bounds.

We now consider an arbitrary system (\ref{1}) and suppose that
it commutes with another system
\refstepcounter{equation}
\label{3}
$$
\begin{array}{ll}
\dot x &= r (x, y),\\
\dot y &= s (x, y).
\end{array}
\eqno{(\ref{3})}
$$
In the language of vector fields, this amounts
to validity of the equalities
\refstepcounter{equation}
\label{4}
$$
\begin{array}{ll}
r p_x + s p_y - p r_x - q r_y &= 0, \\
r q_x + s q_y - p s_x - q s_y &= 0.
\end{array}
\eqno{(\ref{4})}
$$
Then  the above argument clearly remains valid
in each domain where systems (\ref{1}) and (\ref{3})
are transversal. In particular, if such domain surrounds
a center of system (\ref{1}) then the center is isochronous.
It was Sabatini \cite{7} who noticed this fact.
He also proved the converse assertion, namely:
in a punctured neighborhood of an isochronous center
there always exists a system transversal to the original system
and commuting with it; moreover, this system has
the same smoothness as the original system.
This gives a new tool for proving isochronicity of various
systems in a uniform way.

It is worth noting that conditions (\ref{4}) amount merely to
local commutation of the flows $ U_{\tau} $ and $ V_{\sigma} $
corresponding to systems (\ref{1}) and (\ref{3}), even if
the corresponding vector fields are defined on the whole plane.
Surely, sometimes we can extend equality
between the values $ V_{\sigma} U_{\tau} $ and
$ U_{\tau} V_{\sigma} $ at a given point to
rather large $ \tau $ and $ \sigma $, for example, when
the products $ V_{\sigma'} U_{\tau'} $
and $ U_{\tau'} V_{\sigma'} $ make sense at this point
for all pairs $ \sigma', \tau' $ in the rectangle
$ 0 \le \sigma' \le \sigma$,  $0 \le \tau' \le \tau$.
However, we cannot guarantee the equality
$ V_{\sigma} U_{\tau} =
U_{\tau} V_{\sigma} $ on the only reason that
both sides of the equality are defined at the point in question.
Neglect of this circumstance might lead to erroneous conclusions.
For example, by an argument similar to that above, one might
conclude that two isochronous centers whose basins
are joined by an orbit of a commuting system
have the same period. However, this is not always the case.
We exhibit one example.

A system $ \dot z = f(z)$ on the plane,
with $ f(z) $ a holomorphic function of a complex variable $ z $,
has a center at a point $ z_0 $ if and only if $ f(z_0) = 0 $ and
the derivative $ f'(z_0) $ is a nonzero purely imaginary number.
This result is also due to Lukashevich \cite{5}, but
essentially the same result was recently reproved in \cite{8}.
Thus, if, say, $f(z) = i z (1 - z^2) $ then our system
has three centers. They all lie on the real axis, separating
it into four parts each of which, as we easily see,
serves as an orbit
of the orthogonal system $ \dot z = - i f(z)$.
Two of them join the point $ (0, 0) $ to the points
$ (- 1, 0) $ and $ (1, 0) $. Moreover, the extreme centers
have period $ \pi $, whereas the period of the middle center
equals $ 2 \pi $.

We now study the mechanism of this phenomenon in more detail.
The essential elements of the phase portrait of the system
under consideration are displayed in Fig.~1.
Two separatrices infinite to both sides separate the plane
into three domains. These domains
are entirely covered by cycles.
The whole picture is symmetric about the coordinate axes.

\medskip
\centerline{\hfill {\epsfbox{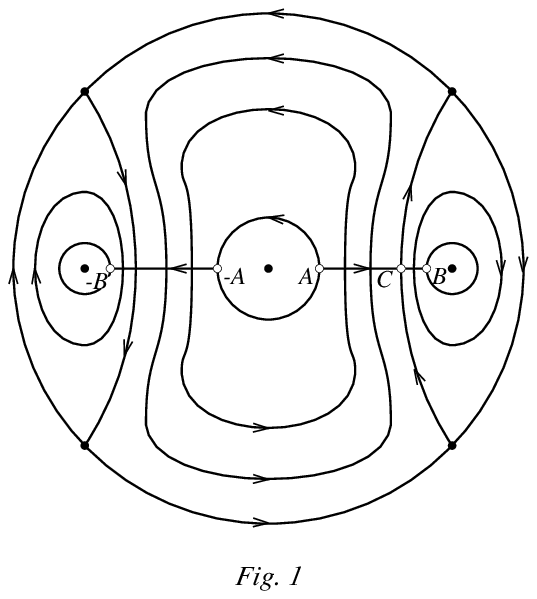}}\quad \quad \quad
{\epsfbox{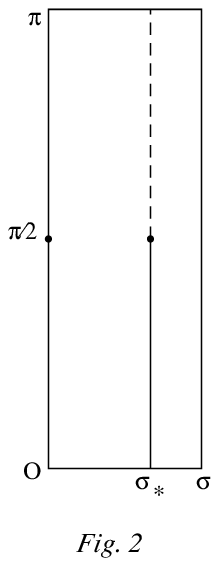}}\hfill}
\medskip
In the orthogonal system, consider the motion
starting at $ t = 0 $ at the point $ A $ indicated in the picture
and finishing at some time $ t = \sigma $ at the point $ B $
lying at the intersection of one of the orbits
around the right center with the abscissa axis.
If we use the former notations $ U $ and $ V $
for the original and orthogonal flows and
take symmetry of the system and the above-indicated values
of periods into account, then we come
to the following expressive relations:
$$
U_{\pi} V_{\sigma} (A) = U_{\pi} (B) = B, \quad
V_{\sigma} U_{\pi} (A) = V_{\sigma} (- A) = - B.
$$

Certainly, the culprits are the separatrices separating
the basins of the centers: their blame consists in
their going to infinity too fast, so that they create obstacles
to extending commutation of the flows to the needed range.
Indeed, we can easily calculate that these separatrices
represent two branches of the hyperbola $ 2 x^2 - 2 y^2 = 1 $
and the motion along them is described by the explicit formulas
$$
x = \pm \frac{\cos(t/2)}{\sqrt{2 \cos t}}, \quad
y = \pm \frac{\sin(t/2)}{\sqrt{2 \cos t}}.
$$
The duration of traversing each of these orbits
embeds in the interval $ - \pi/2 < t < + \pi/2 $ of length~$ \pi $.
Thus, if $ \sigma_{*} $ is the moment when the point
$ C = V_{\sigma_{*}} (A) $ occurs on the right separatrix
then the orbit $ U_t (C) $ terminates by time $ t = \pi/2 $.
Now, disregarding the fact that the other orbits of the flow $ U $
which start at the points of the segment $ A B $
live for the infinite time,
there is no ``commuting''
rectangle which could equate the values of $ V_{\sigma} U_{\pi} $
and $ U_{\pi} V_{\sigma} $ at~$ A $ (Fig.~2).

Seizing an opportunity, we indicate one more curious property
of the periods of polynomial analytical systems.

\begin{theorem}
If a polynomial system of degree $ n \ge 2 $
satisfies the Cauchy--Riemann conditions and has $ n $
singular points such that $n - 1 $ of them are centers, then
the remaining point must be a~center too. Furthermore,
among all singular points there necessarily exist those
around which the rotation is clockwise as well as those
around which the rotation is counterclockwise.
Finally, the sums of periods for the various types of centers
are the same.
\end{theorem}

{\bf Proof.}
If among the roots
$ z_1, \dots , z_n $ of a polynomial $ f(z) $ of degree $ n \ge 2 $
there are $ n - 1 $ simple roots then
it is clear that the remaining root is simple as well.
Furthermore, the following remarkable identity is easy to verify:
$$
\frac{1}{f'(z_1)} + \dots + \frac{1}{f'(z_n)} = 0.
$$
If all points $ z_k $ with $ k = 1, \dots , n - 1 $
are centers of the system $ \dot z = f(z) $
then, as was mentioned above, the numbers
$ f'(z_k) $ are purely imaginary.
But then the preceding identity implies that
the derivative $ f'(z_n) $ must be purely imaginary as well.
Since this derivative differs from zero, the point $ z_n $
must be a~center. It remains to observe that
the absolute value of the real number $ 2 \pi i / f'(z_k) $
equals the period of cycles around the (isochronous)
center $ z_k $ and the sign of this number determines
the direction of the circuit.

\begin{center}
{\large \bf Sabatini's Problem and Devlin's Example}
\end{center}
We now discuss the application of the above
approach to studying isochronous centers of polynomial
dynamical systems. The following problem was stated in
the already-quoted article~\cite{7}:

\medskip
{\bf Problem.} {\sl Is it true that a polynomial system,
having an isochronous center, commutes with some transversal
polynomial system}?

\medskip
For every quadratic system with an isochronous center,
in \cite{9} the author of this question constructed
a commuting polynomial system
whose linear part is orthogonal to the original linear system.
Using {\it Mathematica} \cite{10}, we
examined  all isochronous
polynomial systems described in the articles
we know. It is surprising, but for the majority of them
the answer to the above question turned out to be affirmative.
We outline the results of our inspection.

{\bf Example 1.}
The article \cite{2} contains a full list of
all (up to homothety) cubic systems symmetric
with respect to an isochronous center.
It is as follows:
$$
\begin{array}{llll}
\dot x &= -y -3 x^2 y + y^3,\quad \quad \quad &\dot x &= -y + x^2 y,   \\
\vspace{6pt}
\dot y &=  x + x^3 - 3 x y^2,\quad \quad \quad &\dot y &=  x + x y^2,  \\
\dot x &= -y +3 x^2 y,\quad \quad \quad &\dot x &= -y -3 x^2 y,   \\
\dot y &= x -2 x^3 +9 x y^2,\quad \quad \quad &\dot y &= x + 2 x^3 - 9 x y^2.
\end{array}
$$

Below we indicate the corresponding polynomial systems,
obviously transversal to the original systems at least near zero,
which commute with them as one can check by straightforward
verification of identities (\ref{4}):
$$
\begin{array}{llll}
\dot x &= x+x^3-3x y^2,
\quad \quad \quad
&\dot x &= x - x^3, \\
\vspace{6pt}
\dot y &= y+3x^2 y -y^3,
\quad \quad \quad
&\dot y &= y-x^2 y, \\
\dot x &= x-5x^3 + 6 x^5,
\quad \quad \quad
&\dot x &= x + 5 x^3 +6x^5,\\
\dot y &= y - 9 x^2 y +18 x^4 y,
&\dot y &= y +9x^2 y +18 x^4 y.
\end{array}
$$

{\bf Example 2.}
As is well known \cite{2}, among Kukles's systems
we have exactly one system (up to a ~similarity transformation)
with an isochronous center:
$$
\begin{array}{ll}
\dot x &=- y,\\
\dot y &= x +3 x y + x^3.
\end{array}
$$
For it, there exists a transversal commuting polynomial system:
$$
\begin{array}{ll}
\dot x &= x + x y + x^3,\\
\dot y &= y -x^2 + y^2 - x^4.
\end{array}
$$

As computer analysis shows, the above isochronous system
remains unique even in the full class of Kukles's systems
when the second equation is specified by an arbitrary
polynomial of third degree.
Some relevant remarks can be found, for instance, in \cite{8}.

It is curious to notice that each of the isochronous systems
of the shape
$$
\begin{array}{ll}
\dot x &= - y + p_2  (x, y),\\
\dot y &= \phantom{-} x + q_2  (x, y),
\end{array}
$$
where $ p_2 $ and $ q_2 $ are homogeneous quadratic polynomials
written in the form they are presented in Table~1 of~\cite{2},
in the variables $ X = x$,  $Y = y - p_2 (x, y) $ coincides,
up to a similarly transformation, with the above-written Kukles's system.

Truly ubiquitous, this system plays a magical role in the problem
under study. We encounter it once more in our article, moreover,
at a place where we did not expect to see it at all.

{\bf Example 3.}
Studying the cubic Kolmogorov systems, the authors of \cite{11}
arrived at the following system with an isochronous center at the origin:
$$
\begin{array}{ll}
\dot x &= -y + 2x y - a x^2 y,\\
\dot y &= x - x^2 + y^2  - a x y^2.
\end{array}
$$
This system commutes, for example, with the polynomial system
$$
\begin{array}{ll}
\dot x &= x - x^2 + y^2 - a x y^2,\\
\dot y &= y  - 2 x y - a y^3.
\end{array}
$$

It is noteworthy that, in all three examples of systems, their
direction fields
are symmetric about the abscissa axis, which makes it obvious that
the origin is a center of each of them.
Having constructed transversal commuting systems, we obtain
one more proof for isochronicity of these centers.

{\bf Example 4.}
A new family of isochronous cubic systems was found in \cite{3}.
Slightly changing notations, we can write it as follows:
$$
\begin{array}{ll}
\dot x &= -y - 2 x y + x^2  (\alpha + \alpha x + \beta y),\\
\dot y &= x + x^2 - y^2 + xy (\alpha + \alpha x + \beta y).
\end{array}
$$

For arbitrary real $ \alpha $ and $ \beta $, the origin is a center
of the system, which is easy to check by using the first integrals
indicated in \cite{3}.
Incidentally, everything happened just the other way about:
the authors started with choosing an integral and then wrote a system for it.
 They studied the result by the Darboux linearization method
and proved isochronicity of the center. Another proof
of isochronicity appears on observing that the system in question
commutes with the following:
$$
\begin{array}{ll}
\dot x &= x + x^2 - y^2 + xy (\alpha + \alpha x + \beta y),\\
\dot y &= y + 2 x y + y^2  (\alpha + \alpha x + \beta y).
\end{array}
$$

As we can see, this example contains the previous:
it suffices to take $ \alpha = 0$ and  $\beta = - a $
and change the sign of the variables $ x $ and $ y$.

{\bf Example 5.}
A remarkable class of isochronous hamiltonian systems
relates to area-preserving transformations \cite{12}.
Suppose that functions $ u = u(x, y) $ and $ v = v(x, y) $
satisfy the identity
$$ u_x v_y - u_y v_x \equiv 1.$$
In that event, if $ u = v = 0$, say, at $ x = y = 0$
then the origin is an isochronous center of the system
$$
\begin{array}{ll}
\dot x &= - u u_y - v v_y, \\
\dot y &= \phantom{+} u u_x + v v_x.
\end{array}
$$
This is obvious, since in the coordinates $ u$, $v$ the system
takes a standard form of the harmonic oscillator: $ \dot u = -v$,  $\dot v = u$.
However, for us it is important to indicate an explicit
form of a commuting system which is polynomial
in case $ u $ and $ v $ are polynomials:
$$
\begin{array}{ll}
\dot x &= \phantom{+} u v_y - v u_y, \\
\dot y &= - u v_x + v u_x.
\end{array}
$$
The fact, that this system really commutes with
and is transversal to the original system,
can be verified by straightforward
calculations, but it is simpler to observe that in the coordinates
$ u$, $v $ the system is written as the system $ \dot u = u$,  $\dot v = v$
corresponding to the canonical Euler field $ (u, v)$ obviously
orthogonal to the field $ (-v, u) $ of the harmonic oscillator
and commuting with it.

It is easy to construct as many isochronous hamiltonian systems
as we wish if we observe that the Jacobian of the mapping
$ (u, v) $ is preserved under the replacement of $ u $
by $ u + \varphi(v) $ or $ v $ by $ v + \psi(u)$,
where $ \varphi $ and $ \psi $ are arbitrary smooth functions
of one variable.

We thus see that the hope of a~positive solution to Sabatini's
problem had certain grounds. However, the existence of a commuting
system, even in a small neighborhood of a singular point,
imposes strong requirements on the global behavior
of the system as soon as we talk about the polynomial vector fields.
It is quite possible that these requirements  are met by
all cubic isochronous systems,
but they are already incompatible with the degree of freedom
of systems of forth degree.
One of the circumstances in support of what was said
is reflected in the following simple theorem slightly developing
the arguments of the beginning of the article.

\begin{theorem}
If two polynomial systems commute in some domain and are not collinear
in at least one point of the plane then every center
of each of the systems is isochronous.
\end{theorem}

{\bf Proof.}
The points where a polynomial of two variables
vanishes either cover the whole plane or constitute
a closed nowhere dense set. Therefore, first, our systems
commute on all of the plane and, second, the basin
of every center of each of them includes a cycle
that the other flow, commuting with the given system,
intersects  somewhere transversally.
But then all familiar arguments
apply and demonstrate that all nearby cycles
have the same period. In this case the period,
depending analytically on the orbit, has a constant value
on the whole basin.

Thus, the natural question arises:
Can two centers  of a~polynomial system coexist if one of them is isochronous
and the other is not?
In view of the preceding theorem, a positive answer to this question
would be simultaneously a negative answer to Sabatini's problem.
We found an appropriate example in the article \cite{13},
wherein the system
$$
\begin{array}{ll}
\dot x &=- y - x^4 + 4   x^2 y^2 + y^4,\\
\dot y &= x - 4   x^3 y
\end{array}
$$
is indicated which has one isochronous
and two nonisochronous centers.

When our article was in print, we received from Sabatini
a large, still unpublished survey of isochronous
centers \cite{14} which contains many of the commuting
polynomial systems we consider, as well as a whole series
of new examples.
Interestingly, this survey also points to Devlin's example \cite{13},
moreover, in connection with the same question.
\newpage

\begin{center}
{\large \bf Homogeneous Perturbations of the Harmonic Oscillator}
\end{center}
Besides the above, there are other reasons for which
an isochronous polynomial system may fail to commute
with any nonproportional polynomial system.
Each of the following systems can serve as an example:
\refstepcounter{equation}
\label{5}
$$
\begin{array}{ll}
\dot x &= - y,\\
\dot y &= x + x^{2 m - 1} y (x^2 + y^2),
\end{array}
\eqno{(\ref{5})}
$$
where $ m $ is an arbitrary natural number.
Here we have only one singular point and therefore
the preceding arguments fail.

We have found and investigated this system
before we learned about Sabatini's problem.
Now, when the question has received an answer, we
might confine ourselves merely to mentioning our example.
However, since polynomial vector fields of degree higher than tree
are not studied practically  and any information could be
useful, in this section we conduct a small investigation
of system (\ref{5}) which is curious in many respects.

\begin{theorem}
The origin is an isochronous center of system $(\ref{5})$.
\end{theorem}

{\bf Proof.}
Symmetry of the direction field of $(ref{5})$ about the ordinate axis implies that
the origin is a center of the system. We write the system
in polar coordinates:
$$
\dot \varrho = \varrho^{2m+2} \cos^{2m-1} \vartheta \sin^2 \vartheta,
\quad
\dot \vartheta = 1 + \varrho^{2m+1} \cos^{2m} \vartheta \sin \vartheta.
$$
Differentiating these equations with respect to $ t $,
we can easily verify that the function
$$
t - \vartheta - \frac{\varrho^{2m+1} \cos^{2m+1} \vartheta}
{2 m + 1}
$$
remains constant along each orbit. If the orbit is closed and,
in consequence, surrounds the origin, then during one circuit
the angle $ \vartheta $ receives the increment $ 2 \pi$
and the fractional summand of the above expression
returns to the initial value.
Therefore, the duration of the circuit equals $ 2 \pi$.

\begin{theorem}
The polynomial systems able to commute with system $(\ref{5})$
may differ from $(\ref{5})$ only by a constant factor.
\end{theorem}

{\bf Proof.}
Suppose that the field $ (p, q) $ of system (5)
commutes with a polynomial vector field $ (r, s)$.
When $ p = - y$, as we have in our case, the first
expression in (4) takes the shape

\refstepcounter{equation}
\label{6}
$$
s = y r_x - q r_y.
\eqno{(\ref{6})}
$$
Now, if we rewrite the second relation in (\ref{4})
by using (\ref{6}) and next compare summands
from the viewpoint of their degree in $ y$, not forgetting
that the degree of the polynomial $ q $ is exactly three,
then we find that the degree $ r $ cannot exceed one, so that
$$ r = r_0(x) + r_1(x) y.$$
The further, rather elementary examination leads to
the equalities $ r_0 = 0 $  and $ r_1 = - c$, where $ c $
is a~ constant. Hence, $ r = -c y $ and so $s = c q$.
Thus, the field $ (r, s) $ is collinear to the field $ (p, q) $
with a~ constant coefficient of proportionality.

\medskip
\centerline{\hfill {\epsfbox{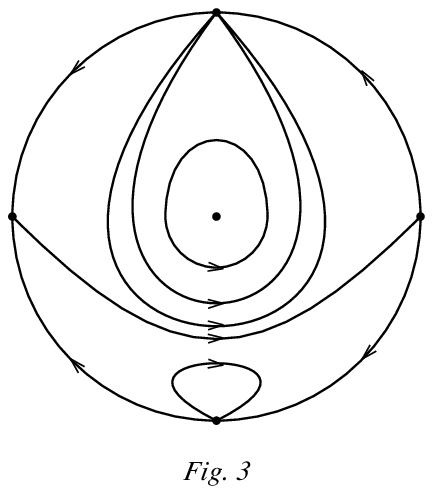}} \hfill}
\medskip
In Fig. 3 we depict the scheme of the phase portrait of system (\ref{5})
which is obtained by the classical methods of the qualitative
theory of planar dynamical systems.
The above example of an isochronous system
resulted from studying the one-dimensional nonconservative
Newton equations with homogeneous nonlinearities
added to the equation of the standard harmonic oscillator.
On the phase plane, these equations are expressed by the system

\refstepcounter{equation}
\label{7}
$$
\begin{array}{ll}
\dot x &= - y,\\
\dot y &= x + Q_n(x, y),
\end{array}
\eqno{(\ref{7})}
$$
where $ Q_n(x, y) $ is a homogeneous polynomial of degree $ n$.
Examining the known full lists \cite{1,2}
of isochronous systems involving only quadratic
or only cubic nonlinearities, we can observe that among
them there are no Newton systems of the form (\ref{7})
with $ n = 2 $ and $ n = 3$. Inspecting the necessary
conditions for isochronicity of a~center by means of {\it Mathematica},
we verified that, up to homothety,
there is exactly one isochronous system for $ n = 4 $:
it is system (\ref{5}) with the exponent $ m = 1$,  and if $ n = 5$
then the origin never can be an isochronous center
of system (\ref{7}). By this we once more confirmed the results
that were obtained thirty years ago in \cite{15}.

In the more difficult case of $ n = 6 $,
we succeeded only in examining systems
that are symmetric about one of the coordinate axes.
We found out that, with the previous stipulations,
the only system among them which has an isochronous center at zero
is system (\ref{5}) with $ m = 2 $.
Finally, if $ n = 7 $ then system (\ref{7}) cannot have an isochronous
center at zero.

The above consideration and some other reasons allow
us to pose the following two questions.

\medskip
{\bf Problem 1.} {\sl
Is it true that a system of the form $(\ref{7})$
with nonlinearities of degree higher than two
has a center at the origin if and only if its direction field
is symmetric about one of the coordinate axes}?

\medskip
{\bf Problem 2.} {\sl Is it true that the origin
is an isochronous center of system $(\ref{7})$ only if the system
has an even degree and is reduced to system $(\ref{5})$ by
a similarity transformation?}

\smallskip
Surely, in near future
the progress in computers will
advance the study of the class of systems
in question far beyond the point we have reached.
However, in this way we could hope to receive
only a negative answer to the posed questions.
If the answers are affirmative then the attempts at their
justification will probably lead to new profound ideas.
These questions are also closely connected with
the other problems to be dealt with below.

\begin{center}
{\large \bf Polynomial Newton Systems}
\end{center}
The circumstances mentioned in the proof of the last theorem
gives us the thought of trying to describe all polynomial Newton systems
for which there are transversal commuting polynomial fields.
The results turned out to be unexpected for us.
The rest of the article is devoted to exposition of these results.

Consider a general polynomial Newton system, expanding its force function
in the powers of the variable having the meaning of velocity:
\refstepcounter{equation}
\label{8}

$$
\begin{array}{ll}
\dot x &= - y,\\
\dot y &= q_0(x) + q_1 (x) y + \dots + q_n(x) y^n.
\end{array}
\eqno{(\ref{8})}
$$
The polynomials $ q_0, \dots , q_n $ in ~$ x $
can have arbitrary degrees, but we assume that $ q_0 $ and $ q_1 $
have no constant terms and that the derivative $ q_0 $
at zero equals one. This means that the linear part of the system
is the canonical oscillator.

As is well known, if a polynomial Newton system (\ref{8})
is conservative, i.e., all $ q_k $ with $ k \ge 1 $
vanish; then the oscillations described by the system
are isochronous only in the trivial case of $ q_0 (x) \equiv x $.
This fact was first established in \cite{16}, afterwards it
was indicated in \cite{1}, and after many years it appeared
in \cite{17}, where it received a new proof.
Thus, in the conservative case, there arises no question
of existence of a commuting system.
Therefore, henceforth we assume that $ n \ge 1 $
and that the last polynomial $ q_n $
does not vanish identically.

\begin{theorem}
If $ n \ge 4 $ then system $(\ref{8})$ cannot commute
with any polynomial system nonproportional to it.
\end{theorem}

{\bf Proof.}
Suppose that a polynomial vector field $ (r, s) $
commutes with the field $ (-y, q) $, where $ q $
is the right-hand side of the second equation in (\ref{8}).
As mentioned, the first relation in (\ref{4})
allows us to express $ s $ by means of (\ref{6}). In this case
the second relation in (\ref{4}) takes the form

\refstepcounter{equation}
\label{9}
$$
q_x (r - y  r_y) + (y q_y - q) r_x + y^2 r_{x x} -
2 y q r_{x y} + q^2 r_{y y} = 0.
\eqno{(\ref{9})}
$$

As in the proof of Theorem 4, this equality implies that
the degree of the polynomial $ r $ in the variable $ y $
cannot exceed one, so that $ r = r_0 (x) + r_1 (x) y$. Thus,
the last summand in (\ref{9}) vanishes and leading is the coefficient
of $ y^{n + 1}$ equal to $ (n - 3) q_n r_1'$.
Since $ n > 3$ and the polynomial $ q_n $ is not equal to zero,
we conclude that $ r_1' = 0 $ or $ r_1 = - c$,  where $ c $
is a constant.

Now, the maximal power of $ y$ in (\ref{9})
equals $ n$ and the corresponding coefficient equals
$$
q_n' r_0 + (n - 1) q_n r_0'.
$$
Its vanishing is possible only if $ (q_n r_0^{n - 1})' = 0 $;
hence, $ r_0 $ too is a constant.
It remains to observe that in these circumstances
equation (\ref{9})
is reduced to the form $ q_x r_0 = 0 $, where $ q_x = 1 $ for $ x = y = 0 $,
which implies that $ r_0 = 0 $. Thus, $r = - c y $,
and so $ s = c q $, as it follows from (\ref{6}).

\medskip
{\bf Problem 3.} {\sl Can a polynomial system of the form $(\ref{8})$
with $ n \ge 4 $ be isochronous}?

\smallskip
Anyway, such systems are absent from the list of the examples
described in the articles known to us.

\begin{theorem}
If $ n = 2 $ then a system of the form $(\ref{8})$
commutes with no polynomial system nonproportional to it.
\end{theorem}

{\bf Proof.}
Here we use the same notations as in the proof of the preceding theorem.
Observe that the highest power with which the variable $ y $
enters in each summand of (\ref{9})  for $ n = 2 $
equals $ m + 2 $, where $ m $ is the degree of $ r $
in $ y $. Hence,
$$
(1 - m) q_2' r_m + (1 - 2 m) q_2 r_m' + r_m'' +
m (m - 1) q_2^2 r_m = 0,
$$
with $ r_m $ the coefficient of $ y^m $
in the expansion of $ r $, representing
a polynomial in $ x $.
Since $ q_2 $ is not identically zero, the last equality implies that
$ m $ is at most one. Certainly, it suffices to settle the case
of $ m = 1 $.
Putting as before $ r = r_0 (x) + r_1 (x) y $,
we note that the equation
is then transformed into $ r_1'' = q_2 r_1' $,
which can hold only for a constant polynomial $ r_1 $.

Collecting the coefficients of $ y^2 $
in (\ref{9}), we now arrive at the relation $ r_0'' = -(q_2 r_0)' $
which means that $ r_0 $ too is a constant.
The further arguments are reduced to repeating what was said
in the proof of the previous theorem: $ r_0 = 0 $
and the vector field $ (r, s) $ is proportional
to the field $ (-y, q) $ of the system under study.

\medskip
{\bf Problem 4.} {\sl
Are there isochronous polynomial systems of the form $(\ref{8})$
for $ n = 2 $}?

\smallskip
We hope that answers to the questions posed in this section
will be related to a rather deep penetration into the
isochronicity nature of polynomial systems.
But already now it should be emphasized that
the very statement of these questions became possible due to
the remarkable idea that had combined concepts
that, as it would seem, have a~little in common:
isochronicity of oscillations and commutation of
vector fields.

\begin{center}
{\large \bf Isochronous Li\'enard Systems}
\end{center}
Of the two remaining classes of polynomial Newton systems,
we first discuss one for which $ n = 1 $.
In this case the force function depends linearly on velocity
and the system, known under the name of a Li\'enard system,
acquires an especially attractive simple form:
\refstepcounter{equation}
\label{10}
$$
\begin{array}{ll}
\dot x &= - y,\\
\dot y &= q_0(x) + q_1 (x) y.
\end{array}
\eqno{(\ref{10})}
$$
We recall that the  linear parts of the systems under study
represent the standard harmonic oscillator.

Here we  digress slightly from the main theme of our investigation,
since in our opinion
the question of commuting fields for systems of the form (\ref{10})
is not very interesting for the following reasons:
In \cite{18} Sabatini studied Li\'enard systems
with odd analytical coefficients $ q_0 $ and $ q_1 $
and obtained an exhaustive solution of
the isochronicity problem in this case.
Namely, he established that under the indicated conditions
the origin is an isochronous center if and only if
$$
q_0 (x) = x + \frac{I^2(x)}{x^3}, \quad {\rm where} \quad
I(x) = \int \limits_0^x s q_1 (s)\,ds.
$$
Furthermore, there are rather solid grounds to believe that,
as far as the polynomial Li\'enard systems are concerned,
there are no other cases of isochronicity, although nobody
can prove this yet (in this connection, see the survey \cite{14}
and the references therein).
Anyway, we are firmly convinced that,
among all polynomial Li\'enard systems, only those
belonging to the Sabatini class can commute
with some transversal polynomial systems; furthermore,
they comprise a small part even in this class.
Clearly, in such situation, the question of commutation
looses actuality, unless its study finds
the end in itself.

We now formulate two problems whose positive solution
would clearly speak in favor of the above-mentioned conjecture.
Assume as before that $q_0 $ and $ q_1 $ are polynomials.
We arrange their primitives:
$$
Q_0 (x) = \int \limits_0^x q_0 (s) \,ds, \quad
Q_1 (x) = \int \limits_0^x q_1 (s)\, ds.
$$

{\bf Problem 5.} {\sl Is it true that
the origin is a center of system $(\ref{10})$ only in the following
two cases: either $ Q_0 $ and $ Q_1 $ are even polynomials
or $ Q_1 $ is a polynomial of $ Q_0 $}?

\smallskip
This problem is cast by a beautiful result by Cherkas
that he obtained in \cite{19}. His theorem, on the one hand,
guarantees the presence of a center under each of the above conditions
and, on the other hand, enables us to reduce the problem to
the following question: is it true that if a polynomial
of one variable is represented as the sum of a power series
in another polynomial which starts with a ~square of the variable
and involves at least one odd power of this variable then
the series terminates somewhere?

\medskip
{\bf Problem 6.} {\sl Is it true that if $ Q_1 $
is a polynomial of $ Q_0 $ then system $(\ref{10}$)
cannot be isochronous}?

\smallskip
As an illustration of the second question, we consider the simplest
class of systems satisfying the indicated condition:
$$
\begin{array}{ll}
\dot x &= - y,\\
\dot y &= f(x) (1 + y),
\end{array}
$$
where $ f(x) $ is a polynomial starting with $ x $.
We prove that the origin is always a center of the system, but
it is never isochronous.

\medskip
\centerline{
{\epsfbox{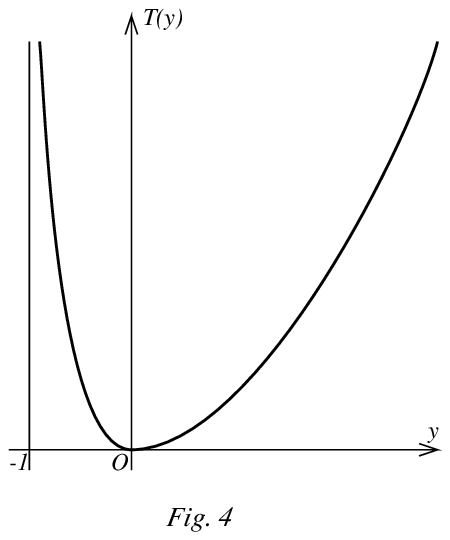}}\quad \quad \quad
{\epsfbox{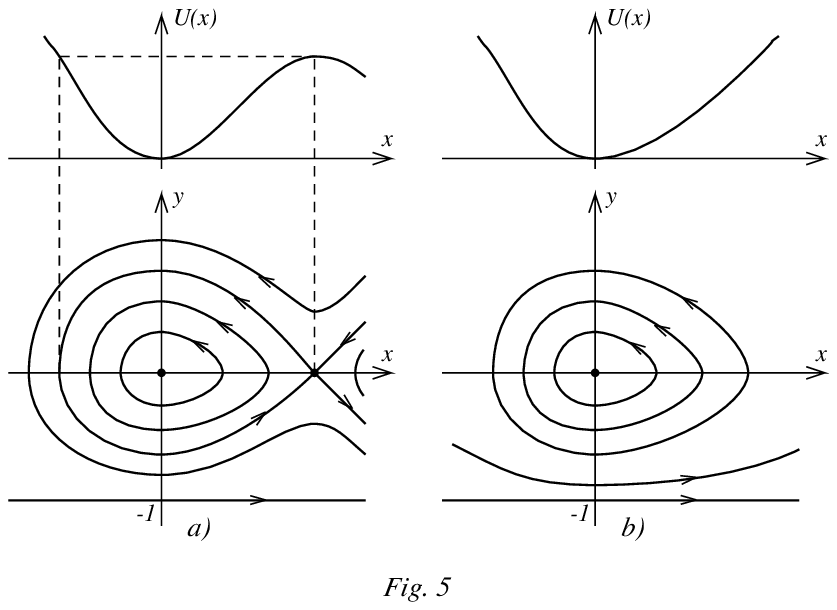}}
}
\medskip
First of all, by setting
$$ U(x) = \int \limits_0^x f(s)\, ds, \quad
T(y) = y - \ln (1 + y),$$
we find that the sum
$$
U(x) + T(y) = \frac{1}{2} (x^2 + y^2) + \dots
$$
serves as a first integral of our system.
Next, we should use the fact that the function
$ T(y) $, defined for $ y > - 1 $,
has two intervals of strict monotonicity (Fig.~4)
on each of which it takes all values from 0 to $ + \infty $,
and examine two possible cases: whether or not the system
has stationary points other than zero.

In the first case, until $ E > 0 $
attains the value of the polynomial $ U $
at the critical point of $U$ nearest to zero
(Fig.~5,\,{\it a}), the equation
$ U(x) + T(y) = E $ determines two functions $ y $ of $ x $
whose graphs form together a cycle of the system.
This means that the singular point nearest to the origin
lies on the boundary of the basin of the center and so
this center cannot be isochronous.
If there are no stationary points other than zero, then
for the same reasons the basin of the center (Fig.~~5,\,{\it b})
adheres immediately to the line $ y = - 1 $
representing an orbit of the system along which the motion
has unit velocity and, consequently, infinite duration,
so that a long time is needed for traversing the nearby
cycles.

The assumption that $ f(x) $ is a polynomial is needed only
for the primitive $ U(x) $ to have two branches going to infinity
in the case of the absence of nonzero roots
and also for the period of a cycle around a center
to depend analytically on the cycle.

\begin{center}
{\large \bf Centralizers of the Abel Systems}
\end{center}
We now consider the polynomial Abel systems, more precisely
systems of the form (\ref{8}) with $ n = 3 $,
slightly changing their record:
\refstepcounter{equation}
\label{11}
$$
\begin{array}{ll}
\dot x &= - y,\\
\dot y &= q_0(x) + 3 q_1 (x) y + q_2 (x) y^2 + q_3 (x) y^3.
\end{array}
\eqno{(\ref{11})}
$$

We do not know all isochronous Abel systems, but
we can give a full description for those of them for which
there are transversal commuting polynomial fields.

\refstepcounter{equation}
\label{12}

\begin{theorem}
The field of system $(\ref{11})$ commutes with some transversal
polynomial vector field if and only if
the following two conditions are satisfied:
$$
\begin{array}{ll}
q_2 q_0 &= 3 q_1^2 - q_0' + 1,\\
q_3 q_0^2 &= q_1^3 - q_1' q_0 +q_1.
\end{array}
\eqno{(\ref{12})}
$$
\end{theorem}

{\bf Proof.}
Let $ (r, s) $ be a polynomial field commuting with the field
$ (-y, q) $ of system (\ref{11}). Then $ s $, as was constantly mentioned,
is expressed by formula (\ref{6}) and $ r $ is found from equation (\ref{9}).
Repeating word for word the arguments presented at the analogous place
in the proof of Theorem 4, incidentally, also devoted to one of
the equations in the Abel class, we conclude that
in this case as well we have $ r =r_0 + r_1 y $, where
$ r_0 $ and $ r_1 $ are polynomials of $ x $.

Elementary calculations of the coefficients of various
powers of $ y $ in (\ref{9}) lead us to four relations, the first two of which
are as follows:
$$
q_0' r_0 - q_0 r_0' =0, \quad
q_1' r_0 -q_0 r_1' = 0.
$$
The general solution of these equation have the form
$$
r_0 =c_0 q_0, \quad r_1 = c_0 q_1 + c_1,
$$
where $ c_0 $ and $ c_1 $ are arbitrary constants.
If $ c_0 = 0 $ then it is easy to see that the commuting field
is proportional to the original field. Therefore,
henceforth we may assume that $ c_0 = 1 $ and $ c_1 = 0 $.

Now, if we write down the third of the above-mentioned four relations,
$$ q_2' r_0 + q_2 r_0' + r_0'' - 6 q_1 r_1' = 0,$$
by using the equalities $ r_0 = q_0 $ and $ r_1 = q_1 $,
we obtain the equation
$$ q_2' q_0 + q_2 q_0' + q_0'' - 6 q_1 q_1' =
(q_2 q_0)' + (q_0')' -3 (q_1^2)' = 0$$
which is easy to solve and which leads to
the formula
$$ q_2 q_0 + q_0' -3 q_1^2 = c_2,$$
where $ c_2 $ is a constant. Factually, $ c_2 = 1 $,
since, as we have agreed, $q_0(0) = q_1(0) = 0 $
and $ q_0'(0) = 1 $. This proves the necessity of
the first condition in (\ref{12}).

To prove the second condition, we use the last relation
which ensues from (\ref{9}):
$$
q_3' r_0 + 2 q_3 r_0' + r_1'' - q_2 r_1' = 0.
$$
Expressing $ r_0 $$, r_1 $ and $ q_2 $ in terms of $ q_0 $
and $ q_1 $, we obtain the differential equation
$$
q_3' q_0^2 + 2 q_3 q_0 q_0' + q_0 q_1'' + q_0' q_1'
- 3 q_1^2 q_1' - q_1' = 0
$$
which integrates without effort:
$$
q_3 q_0^2 + q_0 q_1' - q_1^3 - q_1 = c_3.
$$
Here $ c_3 $ is one more constant. But now it relates
to an expression vanishing at $ x = 0 $, so that $ c_3 = 0 $,
and the necessity of the second condition in (\ref{12})
is established as well.

Sufficiency of (\ref{12}) for existence of a sought commuting
system follows from the obvious conversability
of our logical constructions. Also, it can be easily verified
by inserting the found $ r $ in (\ref{9}).

\refstepcounter{equation}
\label{13}

As we see from the proof of the theorem,
in the case under consideration
there is exactly one sought commuting field $ (r, s) $ which is,
up to a constant factor and a trivial summand, proportional
to the original Abel field $ (-y, q) $:
$$ r = q_0 + q_1 y, \quad s = y r_x - q r_y. \eqno(\ref{13}) $$

Straightforward application of Theorem 7 to construction
of commuting polynomial systems encounters certain difficulties:
if we choose the polynomials $ q_0 $ and $ q_1 $ at random then
the functions $ q_2 $ and $ q_3 $
calculated in accordance with formulas (\ref{12})
turn out as a rule merely rational. Therefore,
it is desirable to have an explicit description for
all polynomial Abel systems whose coefficients are
bound by relations (\ref{12}).
This problem turned out to be solvable.

\refstepcounter{equation}
\label{14}

\begin{theorem}
A polynomial Abel system $(\ref{11})$
commutes with some transversal polynomial system if and only if
it has the following form:
$$
\begin{array}{ll}
\dot x &= - y,\\
\dot y &= (x + a^2 x^3) (1 + h(x) y)^3 +
3 a x y (1 +h(x) y)^2 - h'(x) y^3,
\end{array}
\eqno{(\ref{14})}
$$
where $ a $ is an arbitrary number and $ h(x) $
is an arbitrary polynomial.
\end{theorem}

It is easy to see that the coefficient of $ y^3 $
in (\ref{14}) vanishes only if $ h(x) \equiv 0 $
in which case the system transforms into the isochronous
Kukles system. Even in the simplest case in which $ a = 0 $
and $ h(x) \equiv 1 $, we obtain a new example
of an isochronous Abel system:
$$
\begin{array}{ll}
\dot x &= - y,\\
\dot y &= x (1 + y)^3.
\end{array}
$$
In this case it is obvious that the origin is a center.
Isochronicity of the center is guaranteed
by the membership of the system in the  class (\ref{14}),
although here this fact can be verified by straightforward calculations.

We turn to proving Theorem 8.
In all forthcoming lemmas, the talk is about
the polynomial coefficients of the second equation
of the Abel system (\ref{11}) which satisfy relations (\ref{12})
of Theorem 7.

\refstepcounter{equation}
\label{15}

\begin{lemma}
The polynomial $ q_0 $ has the form
$ q_0(x) = x + A x^3 $, where $ A $ is a constant.
\end{lemma}

{\bf Proof.}
The second condition in (\ref{12}) shows that
the polynomial
$ q_1 \bigl(1 + q_1^2\bigr) $ is divisible by $ q_0 $.
Then it is easy to see that $ q_0 $
admits the decomposition $ q_0 =p_1 p_2 $ into the product
of the polynomials
$$ p_1 (x) = x + O(x^2)  \quad \rm{and} \quad    p_2 (x) = 1 + O(x)$$
which divide $q_1 $ and $ 1 + q_1^2 $ respectively.
Since
$$
q_2 q_0 = 3 q_1^2 - q_0' + 1 = 3 \bigl(q_1^2 + 1\bigr) - q_0' - 2,
$$
the difference $ q_0' - 1 $ is divisible by $ p_1 $
and the sum $ q_0' + 2 $, by $ p_2 $.
Using the equality $ q_0' = p_1' p_2 + p_1 p_2' $,
we arrive at the following important conclusion:
there are polynomials $ h_1 $ and $ h_2 $ such that
$$ p_1' p_2 - 1 = h_1 p_1, \quad
p_1 p_2' + 2 = h_2 p_2. \eqno{(\ref{15})} $$

Before going on further, we exclude one particular case. Namely,
as the first equality in (\ref{15}) shows, if $ p_2 \equiv 1 $
then $ p_1' - 1 = h_1 p_1 $, which is possible only for
$ p_1 \equiv x $; whence, $ q_0 \equiv x $. This is in a~ full agreement
of the claim of the lemma, if we take $ A = 0 $. Thus,
the degree of $ p_2 $ can be assumed positive.
But then, as follows from (\ref{15}), the degrees of $ h_1 $ and $ h_2 $
are less exactly by one than those of $ p_2 $ and $ p_1 $.

\refstepcounter{equation}
\label{16}

We now multiply the first relation in (\ref{15}) by $ h_2 $
and express the resultant product $ h_2 p_2 $ by using
the second relation. Interchanging the words ``one'' and ``two''
in the last proposition, we perform an analogous procedure again.
At a result, we obtain two interesting equalities:
$$
2 p_1' - h_2 = p_1 (h_1 h_2 - p_1' p_2'), \quad
- p_2' + 2  h_1 = p_2 (h_1 h_2 - p_1' p_2').
$$
They are interesting, because they make their left-hand sides
divisible by polynomials of higher degrees and so
yield the formulas
$$ 2 h_1 = p_2', \quad h_2 = 2 p_1'.$$
In this case, both relations of (\ref{15}) become a~sole equality:
$$
2 p_1' p_2 - p_1 p_2' = 2. \eqno (\ref{16}) $$
It is from the last relation that we derive
the necessary conclusions.

First of all, for the mutual compensation of the coefficients
of the highest powers of the variable on the left-hand side of (\ref{16}),
it is necessary that the degree of $ p_2 $, being nonzero, must
be twice as large than the degree of $ p_1 $.
In that event, we can choose a factor $ A $ so that
the degree of $ p_3 = p_2 - A p_1^2 $ be strictly less than
that of $ p_2 $. However, equality (\ref{16}) remains valid
if we substitute $ p_3 $ for $ p_2 $ in it:
$$ 2 p_1' p_3 - p_1 p_3' = 2.$$
These arguments about degrees now lead us to another conclusion:
$ p_3 \equiv c $, where $ c $ is a constant.
Thus, $ p_2 = A p_1^2 + c $, and (16) is extremely simplified:
$ c p_1' = 1 $. Involving the initial data, whence we infer that
$ c = 1 $. Therefore, $ p_1 \equiv x $ and $ p_2 \equiv A x^2 + 1 $,
which completes the proof of the lemma.

\begin{lemma}
There must exist a real number $ a $
and a polynomial $ h (x) $ such that $ A = a^2 $,
where $ A $ is the coefficient mentioned in the preceding lemma
and $ q_1 $ has the form
$$ q_1 (x) = a x + (x + a^2 x^3) h (x).$$
\end{lemma}

{\bf Proof.}
In view of Lemma 1, the first relation in (\ref{12})
looks like
$$ q_2 (x + A x^3) = 3 \bigl(q_1^2 - A x^2\bigr).$$
Were the coefficient $ A $ negative, the left-hand side of the equality
would have nonzero roots, whereas the right-hand side would not.
This proves the possibility of writing $ A = a^2 $.
Then the preceding relation can be written otherwise:
$$ q_2 x (1 + a^2 x^2) = 3 (q_1 - a x) (q_1 + a x).$$
It is clear that one of the nonconstant factors on the right-hand side
must be divisible by the irreducible polynomial $ 1 + a^2 x^2 $.
Since the sign of the number $ a $ was not fixed, we may assume that
this is the first factor.
Moreover, the chosen factor is also divisible by $ x $.
Thus, the function $ h(x) $, necessarily defined as the quotient
of $ q_1 (x) - a x $ by $ x + a^2 x^3 $, is a polynomial.

The next lemma, whose proof is reduced to simple calculations,
is stated only for the reader's convenience.

\begin{lemma}
If $ q_0 $ and $ q_1 $ are expressed in terms of $ a $ and $ h $
as it is said in the previous two lemmas, then equations (12)
is uniquely solvable for $ q_2 $ and $ q_3 $ in the class of polynomials
and their solutions have the form
$$
\begin{array}{ll}
q_2 (x) &= 6 a x h (x) + (x + a^2 x^3) h^2(x), \\
q_3 (x) &= 3 a x h^2(x) + (x + a^2 x^3) h^3(x) - h'(x).
\end{array}
$$
\end{lemma}

To complete the proof of Theorem 8, it now remains only a little:
we have to insert the expressions for coefficients which are indicated
in the three lemmas into the second equation of (\ref{11}) and
verify that the result is system (\ref{14}). Also we should observe that
these coefficients satisfy relations (\ref{12}), because
it is from exactly these relations that the last two of them
were found.

\begin{center}
{\large \bf Conditions for a Center in Smooth Abel Systems}
\end{center}
Describing polynomial Abel systems commuting with transversal polynomial fields,
we supposed nothing about their behavior near the origin in advance.
Meanwhile, if the linear part of a system is a~harmonic
oscillator then the origin can be not only a center but also
a focus. In the latter case, it is quite possible that there exists
a nontrivial commuting field, in which case we also observe some kind
of isochronicity, namely: a point moving along a spiral around a focus
returns to an arbitrarily given radial orbit of the commuting system
by the same time.

Thus, we have to clarify which of the found Abel systems have a center
at the origin. Furthermore, it is interesting to understand the meaning
of the enigmatic form of these systems which is indicated in Theorem 8.
Answers to these questions are given in the following theorem which
holds not only for polynomial systems.

\begin{theorem}
For an arbitrary number $ a $ and an arbitrary smooth function $ h(x) $,
the system of the form (14) in the variables
$$ X = x, \quad Y = \frac{y}{1 + h(x) y} $$
is collinear to the isochronous Kukles system:
$$
\begin{array}{ll}
\lambda \dot X &= - Y,\\
\lambda \dot Y &= \phantom{-} X + 3 a X Y + a^2 X^3,
\end{array}
$$
where $ \lambda = 1 - h(X) Y $. In particular, if $ a = 0 $
then the system is simply reduced to the harmonic oscillator.
Moreover, the origin is always an isochronous center of system $(\ref{14})$.
\end{theorem}

{\bf Proof.}
Verification of reducibility of one system to another
when a change of variables is given explicitly requires nothing
but accurate calculations. Collinearity of a system to the
Kukles system, although with a variable factor,
shows that the origin is a center. Isochronicity
of the center follows from the fact that the field of system (\ref{14})
commutes with the transversal field $ (r, s) $
if we define the latter by formulas (\ref{13}).
This is easy to check by direct calculation of the Poisson bracket or,
equivalently, by verification of equation (\ref{9}).

The just proven theorem can be seen as an assertion
reflecting a curious property of the isochronous Kukles system.
Indeed, if we multiply or divide the field of this system
by an arbitrary function not taking zero values in a neighborhood
of the center then the orbits of the new system are certainly the same
but the law of motion along them changes. However, as we see, if this field
is divided by a function of the form $ 1 - h(X) Y $,
where $ h(X) $ can be arbitrary, then the total time of circuit
around the center remains the same and isochronicity is preserved.
In our opinion, this extraordinary phenomenon which,
judging by all appearances,
should also be encountered in other systems
deserves an independent study.
To state it more precisely, it would be interesting
to clarify the changes of the velocity of the motion along orbits
of an isochronous system which do not influence the period of cycles.

Formulas (\ref{12}) reduce system (\ref{11}) to the form (\ref{14})
only in the case of polynomial vector fields.
If we remove this restriction and consider the Abel systems
with arbitrary smooth coefficients then conditions (\ref{12})
distinguish an essentially wider class of systems.
It turns out that in the general case as well these conditions
not only guarantee existence of a transversal commuting field but also
ensure the presence of a center at the origin.

\begin{theorem}
For arbitrary smooth functions
$ q_0 (x) = x + O(x^2) $ and
$ q_1 (x) = O(x) $, relations $(\ref{12})$
uniquely determine smooth functions $ q_2 (x) $ and $ q_3 (x) $
which are also analytical if so are $ q_0 (x) $ and $ q_1 (x) $.
Moreover, each Abel system (11) with these coefficients
has an isochronous center at the origin.
\end{theorem}

{\bf Proof.}
Since the order of $ q_0 (x) $ in $ x $  is exactly one; to prove
solvability of equations $(\ref{12})$ in the class of smooth or analytical
functions, it suffices to observe that the orders of their right-hand sides
are at least 1 and 2 respectively.

It is curious to notice that systems (\ref{11}) with coefficients
connected by relations (\ref{12}) are contained in the class of integrable
Abel systems which was recently found
by Lukashevich and described by him in \cite{20}.
Following his constructions, we can obtain a first integral of our system
which turns out to be equal to the sum of the squares of two functions
$ X $ and $ Y $ of the variables $ x $ and $ y $.
Omitting the details of intermediate calculations,
we only state the result. Namely, we define these functions
by the formulas
$$
X = x e^{I (x)}, \quad
Y = \frac{x y e^{I (x)}}{q_0 (x) + q_1 (x) y},
$$
where $ I (x) $ stands for the integral
$$
I (x) = \int\limits_0^x \frac{u - q_0 (u)}{u q_0 (u)}\, d u.
$$
Since the integrand clearly has no singularities at zero,
the integral $ I (x) $ represents a smooth function of $ x $.
As it is easy to calculate,
$$ \frac{\partial (X, Y)}{\partial (x, y)} = 1 $$
at the point $ x = y = 0 $, so that the transformation
$ x \to X$,  $y \to Y $ is invertible on some neighborhood
of the origin. Moreover, system (11), with the coefficients
$ q_2 $ and $ q_3 $ expressed in terms of $ q_0 $ and $ q_1 $
by means of (12), in the new variables is reduced to the form
$$
\mu \dot X = - Y, \quad \mu \dot Y = X,
$$
where the variable factor $ \mu $ equals 1 at $ X = Y = 0 $.
Thus, the origin is a center of the system under consideration.

As far as isochronicity is concerned, here it suffices
to repeat what was said in this connection in the proof
of the previous theorem: as it is easy to check,
the field of the system commutes with the field $ (r, s) $  transversal
to it in a neighborhood of the origin and defined by relations (\ref{13}).

However, in this case isochronicity of the center follows also
from other arguments. Namely, the above-mentioned factor $ \mu $
in the new variables $ X $ and $ Y $ has the form $ 1 - f (X)  Y $,
where $ f (X) $ is the ratio $ q_1 (x)/X $
expressed in terms of $ X $, and it is easy to verify
by passage to polar coordinates that such a factor
does not change the period of cycles of the harmonic oscillator.

Closing the article, we return to the polynomial Abel systems
and sum up some conclusions.
First, as Theorem 7 and, especially, Theorem 8 demonstrate,
such systems commute with transversal polynomial fields
only in exceptional cases. Second, the very fact of commutation,
as Theorem 9 or more general Theorem 10 asserts, guarantees
the presence of a center at the origin and this center must be
isochronous. Third, system of the form (\ref{5})
and Theorems 3 and 4 devoted to them recall us that
Theorem 8 describes far from all isochronous polynomial Abel systems.

\medskip
{\bf Problem 7.} {\sl Find all polynomial Abel systems
possessing isochronous centers.}

{\bf Acknowledgement}
We are thankfull to S.Treskov and V.Vershinin for your help during the
preparation of the article.

{\large

Sobolev Institute of Mathematics, 630090, Novosibirsk, Russia.

E-mail: volok@math.nsc.ru
}

\end{document}